\documentstyle[11pt]{article}
\begin{document}
\title{ {Elliptic Curves of Twin-Primes Over Gauss Field and \\
Diophantine Equations }
     \thanks{MR(1991) Subject Classification:
      11R58; 11R27; 14H05; 11A55}
         \thanks{Project Supported by the NNSFC (No. 19771052)}
         }
\author{\mbox{}\vspace{0.3cm}
   {{\Large Q}IU {\Large D}erong 求得容 $(\hskip 1.6cm)\;\;
 and\;\; $ {\Large Z}HANG {\Large X}ianke }(\hskip 1.6cm ) \\
 \small{Tsinghua University,
     Department of Mathematical Sciences,
            Beijing 100084, P. R. China 中国 } }
\date{}
\maketitle
\parindent 24pt
\baselineskip 18pt
\parskip 0pt

 Let $p, q$  be twin prime numbers with $q-p=2$ .
  Consider the elliptic curves
   $$\hskip 2cm E=E_\sigma :\hskip 0.4cm \; y^2 = x (x+\sigma p)
(x+\sigma q)\ . \hskip 1.6cm (\sigma =\pm 1) \hskip 2.2cm (1)  $$
  $\ E=E_\sigma $ is also denoted as $E_+$  or $E_-\ $
   when $\sigma = +1$ or $-1$.
Here the Mordell-Weil group  and the rank
 of the elliptic curve $E$
over the Gauss
field $K={{\bf Q}}(\sqrt -1)$ (and over the rational
field ${{\bf Q}}$) will be determined in several cases; and
results on solutions of
related Diophantine equations and simultaneous Pellian
equations will be given. The arithmetic constructs
over ${{\bf Q}}$ of the elliptic curve $E$
 have been studied
in [1],  the Selmer groups are determined,   results on
 Mordell-Weil group, rank,  Shafarevich-Tate group, and
   torsion subgroups are also obtained. Some results on torsion
 subgroups in [2] will be used here to determine $E(K)$.

   \par Similarly to (1), some other special types of elliptic
   curves  were studied
by A. Bremner, J. Cassels, R. Strocker, J. Top,  B. Buhler,
 B.  Gross and D.  Zagier
 (see [3-5]), $\; e.g.,\;  \; y^2=x(x^2+p), \; $
   $\; y^2=(x+p)(x^2+p^2) , \; $  and
   $\; y^2 = 4x^3 - 28x + 25 .\ $
  The last elliptic curve has
   rank 3 and is famous in solving the Gauss conjecture.
\par

 For any number field  $K$,  the Mordell-Weil theorem asserts
that the  $K-$rational points of $E$ form a finitely generated
 abelian group (the Mordell-Weil group):
 $\ E(K) \cong E(K)_{tors} \oplus {{\bf Z}}^r,\ $
 where $\ E(K)_{tors}\ $  is the torsion subgroup,
 $\ r={{\rm rank}}E(K)\ $ is the rank
 (see [6]).
\par \vskip 0.36cm

{\bf Theorem 1}. Let  $E=E_\sigma$ be the elliptic curve in
(1), $\ K={{\bf Q}}(\sqrt -1)$ the Gauss field.

\hskip0.32cm(a) If $p\equiv 5\ (\mbox{mod}\ 8),\ $ then
 $$\mbox{rank}E(K)=0, \hskip 1cm
E(K)\cong {{{\bf Z}}/2{{\bf Z}}}\oplus
{{{\bf Z}}/2{{\bf Z}} }. $$ \par
\hskip0.32cm(b) If $p\equiv 3\ (\mbox{mod}\ 8),\ $
 $q=a^2+b^2,\ $
$(a+\epsilon)^2+(b+\delta)^2=c^2\ $ ($\epsilon , \delta
=\pm 1,\  a, b, c \in {{\bf Z}}),\ $   then
 $$\hskip 1.2 cm \mbox{rank}E(K)=1, \hskip 1cm
E(K)\cong {{{\bf Z}}/2{{\bf Z}}}\oplus
{{{\bf Z}}/2{{\bf Z}} }\oplus {{\bf Z}}.\hskip 0.5cm $$ \par
\hskip0.32cm(c)\hskip 2.6cm $\mbox{rank} E(K) \leq 3. $ \par
\hskip0.32cm(d) If $p, q$ are as in (b), then
 $$\hskip 1.1cm \mbox{rank}E_+ ({{\bf Q}})=1, \hskip 1cm
E_+({{\bf Q}})\cong {{{\bf Z}}/2{{\bf Z}}}\oplus
{{{\bf Z}}/2{{\bf Z}} }\oplus {{\bf Z}}. $$

In [7], Euler considered the Diophantine equation
$$\left\{
  \begin{array}{c}
    X^2+MY^2=S^2 \\
    X^2+NY^2=T^2 \
  \end{array}\right. $$
 for integers $ M\neq N$, and studied the problem of
 classification of the couples $(M, N) $ such that
 the above equation has non-trivial solutions, which
 is the famous``Euler Concordant Form problem".
 E. Bell, T. Ono, K. Ono studied this problem further
 in [8-10]. \quad  For twin prime numbers $p, q$ we consider
 the following two Diophantine  equations similar to the Euler's :
$$\hskip 5cm  \left\{
  \begin{array}{l}
    X^2-pY^2=S^2 \\
    X^2-qY^2=-T^2 \
  \end{array}\right. \hskip 5cm (I)$$
  $$\hskip 5.2cm   \left\{
  \begin{array}{l}
    X^2-pY^2=2S^2 \\
    X^2-qY^2=-2T^2 \
  \end{array}\right. \hskip 4.8cm (II) $$
If $(X, Y, S, T)$ is a primary solution of (I) (or (II))
(i.e., an integer solution with $(X, Y)= 1 $
 and $XY\neq 0  $), then by [1] we know that
 $(-X^2/Y^2,\ XST/Y^3)$ (or $(-X^2/Y^2,\ 2XST/Y^3)$)
  is a ${{\bf Q}}-$rational
 point of elliptic curve $ E_+ $ in (1). This fact leads us
 to the following theorem.
\par \vskip 0.3cm

 {\bf Theorem 2}. If $p\equiv 5\ (\mbox{mod}\ 8),\ $
 then equations (I) and (II) have no primary solution.
 (Furthermore, the only solution of them is
 $(0, 0, 0, 0)$ )
\par \vskip 0.26cm

For elliptic curve $E=E_+$ in (1) and the Gauss field
$K={{\bf Q}}(\sqrt -1)$, define
$$E({{\bf Q}})^- =
\{(x, y)\in E({{\bf Q}})\ |\ x<0 \}, $$
$$ \hskip 0.2cm n_{2K}^- = \frac{1}{2}\ \# (2E(K)\cap E({{\bf Q}})^-), $$
and let $n(I)$ denote the cardinal of the set of
positive primary solutions of equation (I).

\par \vskip 0.3cm
{\bf Theorem 3} (3.1)\ If $\ n_{2K}^-\neq 0, \ $ then equation
 (I) has a primary solution. \par
\hskip 0.4cm (3.2) We have $\ n_{2K}^- \leq n(I) .\ $
In particular,  if $n_{2K}^- =\infty $
 then equation (I) has
infinitely many primary solutions.
\par \vskip 0.2cm

Finally, consider the simultaneous Pellian equations
$$\hskip 5cm   \left\{
\begin{array}{l}
  x^2-py^2=\sigma \\
  z^2-qy^2=\sigma \hskip 2cm (\sigma=\pm 1)
\end{array} \right. \hskip 3cm (III) $$
This kind of equations has been studied in [11-12]
(but $\sigma =1$ and $\{p, q\}=\{m, n\}$ are general integers).
In particular, Rickert in [12]
obtained:\ if $m=2,\ n=3$  then the
 equation  has no nontrivial integer solution
(Nontrivial means $y\neq 0$). We find that equation (III)
has an interesting relation with elliptic curve $E$ in (1):
a nontrivial integer solution $(a, b, c)$ of (III)
gives a ${{\bf Q}}-$rational point $\ (1/b^2, ac/b^3)\ $
of $E$. We obtain a general result on
equation (III). \par \vskip 0.3cm

{\bf Theorem 4}. (a) If $p\equiv 5\ (\mbox{mod}\ 8),\ $
then equation (III)  has no nontrivial
integer solution.\par
(b) If $p\equiv 3,\ 5\ (\mbox{mod}\ 8), $
then equation (III) with $\sigma =-1\ $ has no nontrivial
integer solution.\par

\vskip 1cm
\begin{center}{\bf {\large \bf R}eferences}\end{center}
 \baselineskip 0pt
\parskip 0pt
\begin{description}

\item[[1]] QIU Derong and ZHANG Xianke,
     Mordell-Weil Groups and Selmer Groups of
     Two Types of Elliptic Curves, to appear
\item[[2]] QIU Derong and ZHANG Xianke, Explicit determination
of torsion subgroups of elliptic curves over multi-quadratic
fields, Advances in Math., 28(1999), 475-478.
\item[[3]] A. Bremner and J. W. S. Cassels, On the equation
 $ y^2=x(x^2+p),\ $   Math. Comp. 42(1984), 257-264.
\item[[4]] R. J. Stroeker and J. Top, On the equation
$ y^2=(x+p)(x^2+p^2),\ $ , Rocky Mountain J. of Math.
      24(1994), 1135-1161.
\item[[5]] J. P. Buhler, B. H. Gross and D. B. Zagier, On the
conjecture of Birch and Swinnerton-Dyer for an elliptic curve of rank 3, Math. Comp.
44(1985), 473-481.
\item[[6]] J. H. Silverman, The Arithmetic of Elliptic Curves,
 GTM 106, Springer-Verlag, 1986.
\item[[7]] L. Euler, De binis formulis speciei
$xx+myy$ et $xx+nyy$ inter se concordibus et disconcordibus,
Opera Omnia Series I, 5(1780), 48-60,
Leipzig-Berlin-Z$\ddot{u}$rich, 1944.
\item[[8]] E. T. Bell, The problems of congruent numbers and
concordant forms, Proc. Amer. Acad. Sci., 33(1947), 326-328.
\item[[9]] T. Ono, Variations on a theme of Euler,
Plenum, New York, 1994.
\item[[10]] K. Ono, Euler's condordant forms, Acta Arithmetica,
LXXVIII. 2(1996), 101-123.
\item[[11]] H. P. Schlickewei, The number of subspaces
occurring in the p-adic subspace theorem in Diophantine
approximation, J, Reine Angew. Math., 406(1990), 44-108.
\item[[12]] J. Rickert, Simultaneous rational approximations
and related Diophantine equations, Math. Proc. Cambridge
Philos. Soc., 113(1993), 461-472.

\end{description}
\par  \vskip 0.3cm

{\Large T}SINGHUA {\Large U}NIVERITY \par \vskip 0.136cm {\Large D}EPARTMENT OF {\Large
M}ATHEMATICAL {\Large S}CIENCES \par \vskip 0.136cm {\Large B}EIJING 100084,  P. R.
{\Large C}HINA \par \vskip 0.16cm

E-mail:\hskip 0.3cm  xianke@tsinghua.edu.cn

\end{document}